\documentclass[11pt]{article}
\usepackage{graphicx}
\usepackage{subfigure}
\usepackage{tikz}
\usepackage{tikz-network}
\usepackage{float}
\usetikzlibrary{positioning,petri}
\usepackage{amsfonts}
\usepackage{amsbsy}
\usepackage{amssymb}
\usepackage{amsmath,amsthm}
\usepackage{verbatim}
\usepackage{amsfonts}
\usepackage{parskip}
\usepackage[T1]{fontenc}
\usepackage{multicol}
\usepackage[onehalfspacing]{setspace}
\usepackage{color}
\usepackage[autostyle]{csquotes}
\usepackage{pgf}
\usepackage{mathtools}
\usepackage{hyperref}
\usepackage{thmtools}

\declaretheoremstyle[%
  spaceabove=-6pt,%
  spacebelow=6pt,%
  headfont=\normalfont\itshape,%
  postheadspace=1em,%
  qed=\qedsymbol%
]{mystyle} 
\declaretheorem[name={Proof},style=mystyle,unnumbered,
]{prf}
\hypersetup{
    colorlinks=true,
    linkcolor=blue,
    filecolor=blue,      
    urlcolor=blue,
    citecolor=blue,
    pdftitle={Overleaf Example},
    pdfpagemode=FullScreen,
    }

\newtheorem{thm}{Theorem}
\newtheorem{prop}{Proposition}[section]
\newtheorem{cor}[thm]{Corollary}
\newtheorem{defn}[thm]{Definition}
\newtheorem{obs}[prop]{Observation}
\newtheorem{lem}[thm]{Lemma}

\newtheorem{exm}[thm]{Example}

\newcommand{\Z}{\mathbb{Z}}

\RequirePackage{pgfkeys}
\pgfkeys{/ytableau/options/.is family}

\pgfkeys{/ytableau/options, boxsize/.value required,  boxsize/.code = {%
\pgfkeysalso{nosmalltableaux}%
 \compare@yT{#1}{normal}%
 \ifeq@yT
 \xdef\macro@boxdim@yT{\expandonce@yT\boxdim@normal@yT}%
 \else
 \xdef\macro@boxdim@yT{#1}%
 \fi
 } }
 \pgfkeys{/ytableau/options, textmode/.value forbidden,
 textmode/.code = \set@textmode@yT,
mathmode/.value forbidden,
mathmode/.code = \set@mathmode@yT,
 }
 
\def\set@mathmode@yT{
 \gdef\skipin@yT{$}
 \gdef\skipout@yT{$}
 \def\smallfont@yT{\scriptstyle} } 
 \set@mathmode@yT

\makeatletter
\let\@fnsymbol\@arabic
\makeatother

\title{Realization of zero-divisor graphs of finite commutative rings as threshold graphs}
\author{Rameez Raja$^{1}$ and Samir Ahmad Wagay$^{2}$\footnote{$^{, 2}$Department of Mathematics, National Institute of Technology Srinagar, Jammu and Kashmir, India. Email: rameeznaqash@nitsri.ac.in, samir\_03phd20@nitsri.net,
Corresponding author:$^1$}}

\begin{document}
\maketitle
\begin{abstract} 
  Let $R$ be a finite commutative ring with unity, and let $G = (V, E)$ be a simple graph. The zero-divisor graph, denoted by $\Gamma(R)$ is a simple graph with vertex set as $R$, and two vertices $x, y \in R$ are adjacent in $\Gamma(R)$ if and only if $xy=0$. In \cite{SP}, the authors have studied the Laplacian eigenvalues of the graph $\Gamma(\Z_{n})$ and for distinct proper divisors $d_1, d_2, \cdots, d_k$ of $n$, they defined the sets as, $\mathcal{A}_{d_i} = \{x \in \Z_{n} : (x, n) = d_i\}$, where $(x, n)$ denotes the greatest common divisor of $x$ and $n$. In this paper, we show that the sets $\mathcal{A}_{d_i}$, $1 \leq i \leq k$ are actually orbits of the group action: $Aut(\Gamma(R)) \times R \longrightarrow R$, where $Aut(\Gamma(R))$ denotes the automorphism group of $\Gamma(R)$. Our main objective is to determine new classes of threshold graphs, since these graphs play an important role in several applied areas. For a reduced ring $R$, we prove that $\Gamma(R)$ is a connected \textit{threshold} graph if and only if $R\cong {F}_{q}$ or $R\cong {F}_2 \times {F}_{q}$. We provide classes of \textit{threshold} graphs realized by some classes of local rings. Finally, we characterize all finite commutative rings with unity of which zero-divisor graphs are not \textit{threshold}. 
\end{abstract}


2020 AMS Classification Code: 13A70, 05C25, 05C50.  
\section{Introduction}
The zero-divisor graph of a commutative ring $R$ with unity was introduced by Beck \cite{Bk}. The graph $\Gamma(R)$ is defined to be the graph with vertex set as $R$, and two distinct vertices $x, y\in R$ are adjacent in $\Gamma(R)$ if and only if $xy = 0$ in $R$. The notion of zero-divisor graphs was further studied by Anderson and Livingston in \cite{AdLs}, who restricted  vertices of the graph $\Gamma(R)$ to non-zero zero-divisors of the ring $R$. The main objective of the investigation of associating a graph to $R$ is to study the interplay of graph theoretic properties of $\Gamma(R)$ and ring theoretic properties of $R$. This interplay between different properties of $R$ and $\Gamma(R)$ was thoroughly investigated in \cite{B, PBT, PBS, RR, Rd}. For basic definitions from graph theory, we refer to \cite{SP}, and for ring theory we refer to \cite{K}.

A vertex in a graph $G$ is called \textit{dominating} if it is adjacent to every other vertex of $G$. A graph $\mathbf{G}$ is called a \textit{threshold} graph if it is obtained by the following procedure.

Start with $K_1$, a single vertex, and use any of the following steps, in any order, an arbitrary number of times.\\
(i) Add an isolated vertex.\\
(ii) Add a dominating vertex, that is, add a new vertex and make it adjacent to each existing vertex. The family of threshold graphs represent an important class of simple graphs.
  
The concept of threshold graphs was first established by the authors Chvátal and Hammer in \cite{CH} and Henderson and Zalcstein in \cite{HZ}. There are many different generalizations of threshold graphs. These graphs can be viewed as special cases of some wider classes of graphs like cographs, split graphs and interval graphs.
  
By the definition of Beck's zero-divisor graph $\Gamma(R)$, it appears that the vertex $0$ in $\Gamma(R)$ is adjacent to every other vertex of $\Gamma(R)$. Therefore,  $\Gamma(R)$ is a simple connected graph, that is,  $\Gamma(R)$ is without parallel edges and self loops. The distance between any two vertices in $\Gamma(R)$ is at least $1$, which is the minimum eccentricity of $\Gamma(R)$ and at most $2$,  which is the maximum eccentricity of $\Gamma(R)$. 

A graph $G = (V, E)$ is said to be complete if every pair of distinct vertices of $G$ are adjacent in $G$. Recall that a \textit{clique} of $G$ is a complete subgraph of $G$. A subset $S$ of a vertex set of $G$ is said to be an \textit{independent set} if no two members of $S$ are adjacent in $G$. A graph $G$ is said to be \textit{bipartite} if $V$ can be partitioned into non-empty subsets $V_1$ and $V_2$ such that each edge of $G$ has one end in $V_1$ and other in $V_2$. A bipartite graph is said to be complete if every element of $V_1$ is adjacent to every element of $V_2$. A complete bipartite  is called a \textit{star graph} if either $|V_1| = 1$ or $|V_2| = 1$. If $R\cong {F}_{q}$, where ${F}_q$ is a finite field, then it is easy verify that $\Gamma(R)$ is a star graph 
  
An element $r \in R$ is said to be \textit{nilpotent} if $r^n = 0$, for some $n \in \mathbb{Z}_{>0}$. A ring $R$ is called \textit{reduced} if it has no non-zero nilpotent element and it is called \textit{local} if $R$ has a unique maximal ideal. 
   
This research article is organised as follows. In Section 2, we present some preliminary information related to the zero-divisor graph $\Gamma(\Z_{n})$. For $R\cong \mathbb{Z}_{p^{\alpha}}$, we show that the sets $\mathcal{A}_{d_{i}} = \{x \in \Z_{n} : (x, n) = d_i\}$ defined in \cite{SP} are actually orbits of the group action $Aut(\Gamma(R)) \times R \longrightarrow R$. In Section $3$, we show that for a reduced ring $R$, $\Gamma(R)$ is a connected threshold graph if and only if $R\cong {F}_{q}$ or $R\cong {F}_2 \times {F}_{q}$. We discuss classes of threshold graphs realized by some classes of local rings and finally we characterize all finite commutative rings of which zero-divisor graphs are not threshold.

\section{Preliminaries}
It is always interesting to determine new classes of threshold graphs, since a threshold graph $\mathbf{G}$ with $n$ vertices can be represented by a binary code of the type $b =(b_1b_2 \cdots b_n)=0^{s_1}1^{t_1}0^{s_2}1^{t_2}\cdots 0^{s_k}1^{t_k}$, where $s_i~~and~~t_i$ are positive integers for all $1\leq i\leq k$. In the binary code $b$, $b_i = 0$ if a vertex $v_i$ is added as an isolated vertex in $\mathbf{G}$ and $b_i =  1$ if $v_i$ is added as a dominating vertex in $\mathbf{G}$. These codes realized by some class of graphs (threshold graphs) are known as the creation sequences of these graphs.

We partition the vertex set $V$ of $\mathbf{G}$ as $V_1, V_2, \cdots, V_k$. The equivalent form of an equitable partition matrix $\mathcal{M}=(a_{ij})$ of $\mathbf{G}$ is defined as,
\begin{equation}
\label{eqone}
a_{ij}=\begin{cases}
        |V_{i}| - 1  & \text{ if } i=j \mbox{ and} ~V_i ~\mbox{induces a complete subgraph in} ~\mathbf{G},\\   
       |V_j| & \text{if } \mbox{there exists} ~v_i \in V_i ~\mbox{and}  ~v_j \in V_j ~\mbox{such that} ~v_i ~\mbox{is adjacent to} ~v_j,\\
           0 &\text{otherwise.}
\end{cases}     
\end{equation}\\
The following example illustrates the equitable partition matrix of a threshold graph. 

\begin{figure}[H]
\centering
\includegraphics[scale=.11]{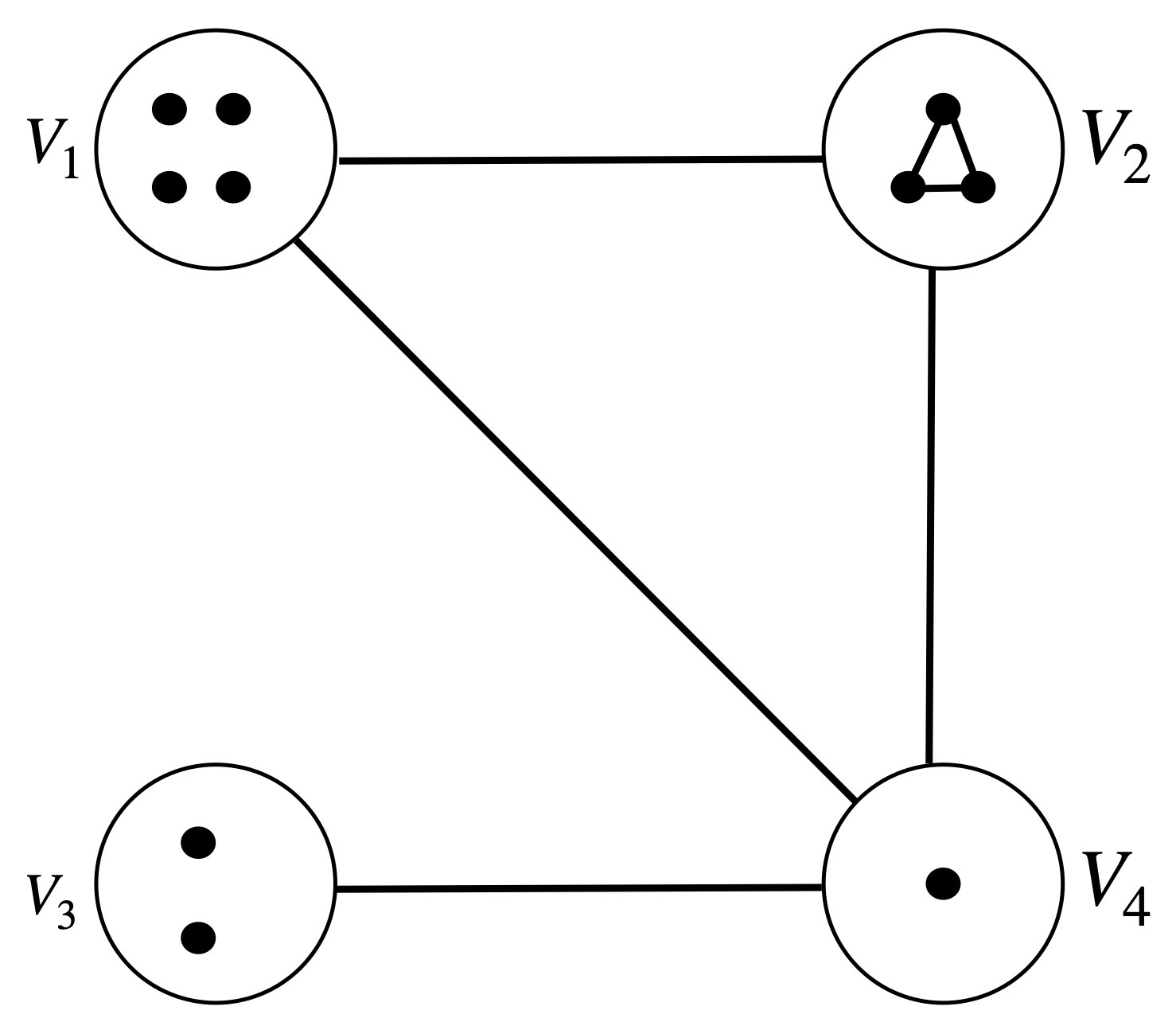}
\caption{Threshold graph (0000111001)}
\label{a}
\end{figure}

Let $\mathbf{G}$ be a threshold graph of order $10$ represented by a binary code $b = (0000111001) = 0^{4} 1^{3} 0^{2} 1$ as shown in the Figure \ref{a}. Then the equitable partition matrix of  $\mathbf{G}$ is,
\begin{equation*}
\mathcal{M} = 
\begin{pmatrix}
0&3&0&1\\
4&2&0&1\\
0&0&0&1\\
4&3&2&0
\end{pmatrix}
\end{equation*}
The multiplicities of eigenvalues $0$ and $-1$ are $4$ and $2$, which can be verified from \cite{RB}. The remaining eigenvalues of $\mathbf{G}$ can be computed from the equitable partition matrix defined in \eqref{eqone}. In fact, the remaining eigenvalues of $\mathbf{G}$ are roots of the characteristic polynomial $x^4-2x^3-21x^2-12x+24$.

Thus constructing some new classes of threshold graphs provides us a way of developing some new binary codes from graphs and a theory which involves spectral analysis related to graphs.

 Recall that the graph obtained from the union of two graphs $G_1$ and $G_2$ by adding new edges from each vertex of $G_1$ to every vertex of $G_2$ is called join of $G_1$ and $G_2$, denoted by $G_1\vee G_2$. The generalised join graph which is also called generalised composition graph in \cite{AJS} is defined as follows.

\begin{defn} Let $G$ be a graph with vertex set $V(G) =\lbrace v_1, v_2, \cdots, v_k\rbrace$, and let $G_1, G_2, \cdots, G_k$ be $k$ pairwise disjoint graphs. The $G-$generalized join graph $G_{\vee} = [G_1, G_2, \cdots, G_k]$ of $G_1, G_2, \cdots , G_k$ is the graph obtained by replacing each vertex $v_i$ of $G$ with the graph $G_i$ and two vertices $G_i$ and $G_j$ are adjacent in $G_{\vee}$ if and only if every vertex $v_i$ of $G_i$ is adjacent to every vertex $v_j$ of $G_j$ in $G$.
\end{defn}
Let $Aut(G)$ denotes the automorphism group of a graph $G$, consider the group action $Aut(G)\\
 ~acting~on ~V(G)$ by some permutation of $Aut(G)$, that is, $Aut(G) \times V(G) \rightarrow V(G)$ is given as, 
 $\sigma(v) = u$, where $\sigma \in Aut(G)$ and $v, u\in V(G)$ are any two vertices of $G$.

The graph $\Gamma(\mathbb{Z}_n)$ is simple and its vertices are all elements of the ring $\mathbb{Z}_n$ and two distinct vertices are adjacent in $\Gamma(\mathbb{Z}_n)$ if and only if there product is zero in $\mathbb{Z}_n$. The canonical decomposition of $n$ is $p_1^{n_1}p_2^{n_2}\cdots p_r^{n_r}$ and total number of divisors of $n$ is $\prod\limits_{i = 1}^{r}(n_i+1)$.

The authors in \cite{SP} defined $\Upsilon_n$ as a simple graph whose vertices are the proper divisors of $n$, and two distinct vertices $d_i$ and $d_j$ are adjacent in $\Upsilon_n$ if and only if  $d_i d_j$ is a multiple of $n$. Thus, $V(\Upsilon_n) = \{d_1, d_2, \cdots, d_r\}$ and 
$E(\Upsilon_n) = \{(d_i, d_j) : d_id_j \equiv 0(mod~n)\}$.
 
Let $n = p^{\alpha}$, where $p$ is a prime number and $\alpha \geq 1$. If $R \cong \Z_{p^{\alpha}}$, then orbits of the group action $Aut(\Gamma(R)) \times R \longrightarrow R$ are given as, $\mathcal{O}_{\alpha,p^i} = \{p^ib (\bmod ~p^{\alpha})\mid b \in \Z, (b,p) =1\}$, $\mathcal{O}_{\alpha,p^i}$ is the $p^i$-th orbit of $R$, where $i \in [0, \alpha-1]$.

For $0\leq i< j \leq \alpha-1$, $p^ib \equiv p^j b' (\bmod ~p^{\alpha}),\text{ where } (b,p) =1 \text{ and } (b',p)=1$.
Therefore, $p^{(j-i)}b' \equiv b (\bmod ~p^{(\alpha-i)}),$ which is a contradiction. Thus, for $i \neq j$, $\mathcal{O}_{\alpha,p^i} \cap  \mathcal{O}_{\alpha,p^j} = \emptyset$.

Any element $a \in \Z_{p^{\alpha}}$ can be expressed as, $$a \equiv p^{\alpha-1}b_1 + p^{\alpha-2} b_2 +\cdots + p b_{\alpha-1 } +b_{\alpha} (\bmod~ p^{\alpha}),$$
where $b_i \in [0, p-1].$ If $a \in \mathcal{O}_{\alpha,1}$,
then $b_{\alpha} \neq 0.$ So, $|\mathcal{O}_{\alpha,1}| = p^{\alpha-1} (p-1) = \phi (p^{\alpha})$. If $a' \in \mathcal{O}_{\alpha,p}$, then for some $a \in \mathcal{O}_{\alpha,1}$, $a' = pa$, that is, $b_{\alpha}\neq 0$, so $|\mathcal{O}_{\alpha,p}| = \phi(\frac{p^{\alpha}}{p})$. Similarly, for $i \in [0,\alpha -1]$, we have $|\mathcal{O}_{\alpha,p^i}| = \phi(\frac{p^{\alpha}}{p^i})$ and  $\mathcal{O}_{\alpha,p^{\alpha}} = \{0\}$. It follows that the size of the sets $\mathcal{O}_{\alpha,p^{i}}=\phi(\frac{p^\alpha}{p^i})$ for $i \in [0,\alpha]$.  Therefore, total number of elements in all the sets $\mathcal{O}_{\alpha,p^{i}}=\sum\limits_{i=0}^{\alpha}\phi(\frac{p^\alpha}{p^i})=|\mathbb{Z}_{p^\alpha}|={p^\alpha}$

Note that $1, p, \cdots, p^{\alpha}$ are the representatives of orbits $\mathcal{O}_{\alpha,1}, \mathcal{O}_{\alpha,p}, \cdots, \mathcal{O}_{\alpha, p^{\alpha}}$ and divisors of the order of $R$. It follows that the vertices of the graph $\Upsilon_n$ are elements of sets $\mathcal{O}_{\alpha,p}, \cdots, \mathcal{O}_{\alpha, p^{\alpha-1}}$. Moreover, by  definition of the graph $\Gamma(R)$, an element $0$ of $R$ is adjacent to all vertices in $\Gamma (R)$, elements of the orbit $\mathcal{O}_{\alpha,1}$ are adjacent to $0$ only in $\Gamma (R)$. Also, elements of the orbit $\mathcal{O}_{\alpha,p}$ are adjacent to $0$ and elements of the orbit $\mathcal{O}_{\alpha,p^{\alpha -1}}$. The elements of the orbit $\mathcal{O}_{\alpha,p^2}$ are adjacent to $0$ and elements of orbits $\mathcal{O}_{\alpha,p^{\alpha -1}}$,  $\mathcal{O}_{\alpha,p^{\alpha -2}}$. This, for $k \geq 1$, elements of the orbit $\mathcal{O}_{\alpha,p^k}$ are adjacent to $0$ and elements of orbits $\mathcal{O}_{\alpha,p^{\alpha -k}}$,  $\mathcal{O}_{\alpha,p^{\alpha -k + 1}}, \cdots, \mathcal{O}_{\alpha,p^{\alpha -1}}$.
 
Thus, the graph $\Upsilon_n$ can be obtained as a subgraph of $\Gamma (R)$ with vertices as elements of orbits $\mathcal{O}_{\alpha,p}, \cdots, \mathcal{O}_{\alpha, p^{\alpha-1}}$, and two vertices $x$ and $y$ are adjacent in $\Upsilon_n$ if and only $xy = 0$.

There is an advantage for knowing the orbits of this group action $Aut(\Gamma(R)) \times R \rightarrow R$, since we get some structural information about some elements of $R$ from $\Gamma(R)$. As a consequence, we do not consider all elements of $R$ to decode the symmetry of $\Gamma(R)$. We explore this information to reveal some interesting spectral properties of  $\Gamma(R)$.

 From the structure of $\Gamma(R)$, the subgraphs induced by elements of orbits  $\bigcup\limits_{i=0}^{k-1} \mathcal{O}_{\alpha,p^{i}}$, are independent parts in $\Gamma(R)$, whereas the subgraphs induced by elements of orbits $\bigcup\limits_{j=0}^{k-1} \mathcal{O}_{\alpha,p^{k
 +j}}$ are cliques in $\Gamma(R)$. Also the orbit $\mathcal{O}_{\alpha,p^{\alpha}}=\lbrace 0 \rbrace$ is connected to every element of $\Gamma(R)$. The following Lemma shows adjacency of vertices in $\Gamma(\mathbb{Z}_{p^{\alpha}})$.
 
 \begin{lem}
 For $i,j \in \lbrace 0, 1, \cdots, \alpha-1 \rbrace$, two vertices $x \in \mathcal{O}_{\alpha,p^{i}}$, $y \in \mathcal{O}_{\alpha,p^{j}}$ of $\Gamma(R)$ are adjacent in  $\Gamma(R)$ if and only if  $i+j \geq \alpha$.
 \begin{prf}
 Since $x\in \mathcal{O}_{\alpha,p^{i}}$ and $y\in \mathcal{O}_{\alpha,p^{j}}$, therefore, $x=p^ib$ and $y=p^jb'$ for some integers $b$ and $b'$ with $(p, b) = (p, b') = 1$. The vertices $x$ and $y$ are adjacent in $\Gamma(R)$ if and only if $p^{\alpha}$ divides $xy$. This implies $x$ and $y$ are adjacent in $\Gamma(R)$ if and only if $i+j \geq \alpha$.
 \end{prf}
 \label{20} 
\end{lem}
\begin{lem}
Let $\Gamma(R)$ be the graph realized by $R \cong \Z_{p^\alpha}$. Then the induced subgraph $\Gamma(\mathcal{O}_{\alpha,p^i})$ of $\Gamma(R)$ on the vertex set $\mathcal{O}_{\alpha,p^i}$ for $i\in\lbrace 0, 1, 2, \cdots, \alpha-1\rbrace$ is either the complete graph $K_{\phi (\frac{p^{\alpha}}{p^{i}})}$ or its complement graph $\overline{K}_{\phi (\frac{p^{\alpha}}{p^{i}})}$. Furthermore, for $\alpha = 2k$ or $2k-1$, $\Gamma(\mathcal{O}_{\alpha,p^i})= K_{\phi (\frac{p^{\alpha}}{p^{i}})}$ if and only if  $i \geq k$ and $\Gamma(\mathcal{O}_{\alpha,p^{\alpha}}) = K_{1}$.
 \end{lem}
\begin{prf}
 The subgraph $\Gamma(\mathcal{O}_{\alpha,p^i})$ has $\phi(\frac{p^\alpha}{p^i})$ vertices for $i \in [0,\alpha]$. By taking $i=j$ in Lemma \eqref{20}, we get the result.
\end{prf}
\begin{obs}
If ~$\Gamma(\mathcal{O}_{\alpha,p^i})$ is the induced subgraph of $\Gamma(R)$ on vertices $\mathcal{O}_{\alpha,p^i}$ for $i \in \lbrace 0, 1, 2, \cdots, \alpha \rbrace $, then $\Gamma(R)_{\vee} = [\Gamma(\mathcal{O}_{\alpha, 1}), \Gamma(\mathcal{O}_{\alpha,p}), \cdots,  \Gamma(\mathcal{O}_{\alpha,p^{\alpha}})]$.
\label{4}
\end{obs}
\begin{prf}
Proof follows from the adjacency relations in $\Gamma(R)$, since orbits are adjacent in $\Gamma(R)_{\vee}$ if and only if their representatives are adjacent in $\Gamma(R)$.  
\end{prf}

We conclude this section  with the following example which illustrates the zero-divisor graph realized by $\Z_{3^3}$ with its orbits.

\begin{figure}[H]
\centering
\includegraphics[scale=.2]{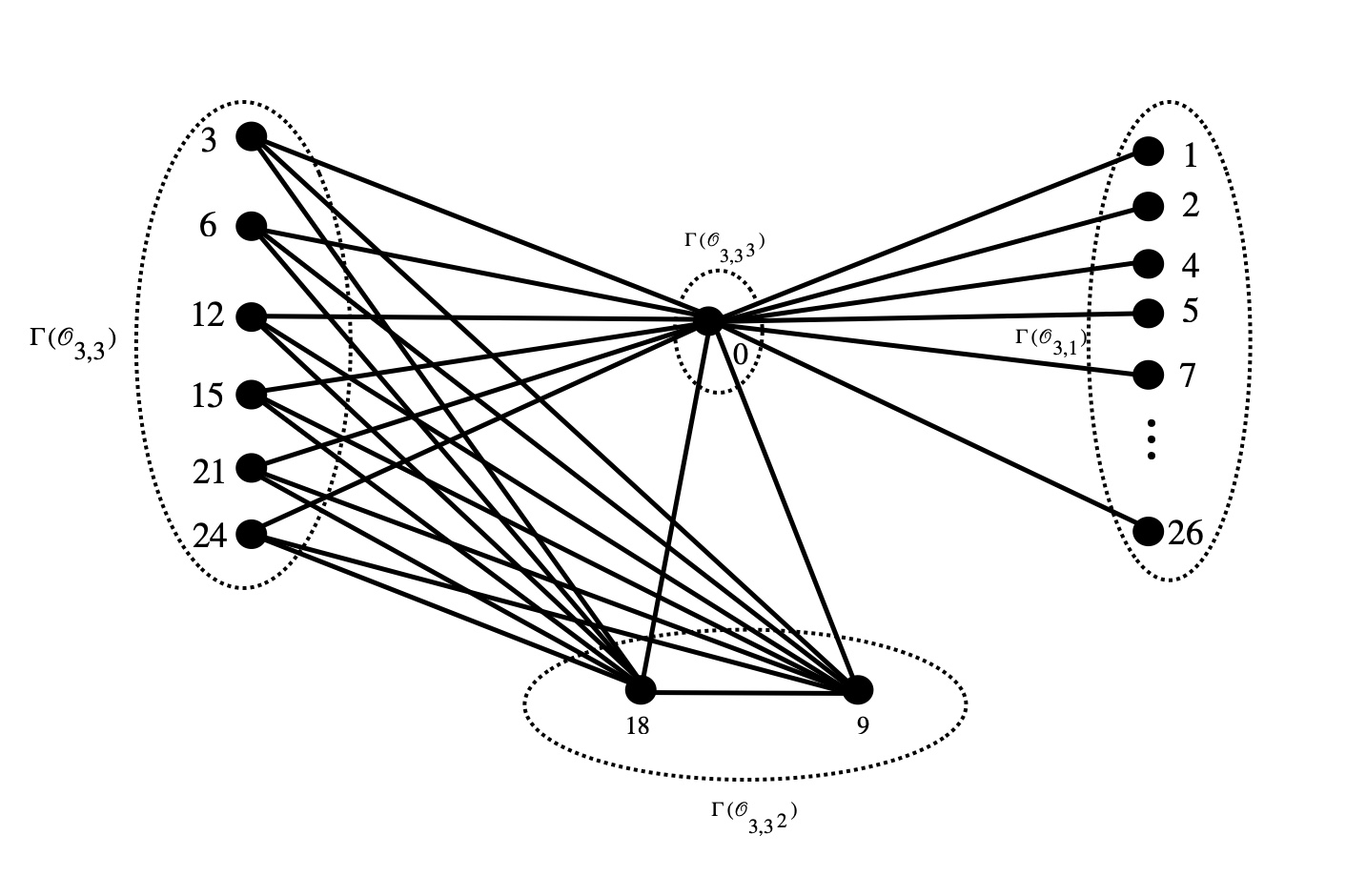}
\caption{Generalised join graph $\Gamma(\mathbb{Z}_{3^3})_{\vee} $}
\label{b}
\end{figure}

\begin{exm}Let $R=\mathbb{Z}_{3^{3}}$. Consider the group action $Aut(\Gamma(\mathbb{Z}_{3^{3}})) \times \mathbb{Z}_{3^{3}} \rightarrow \mathbb{Z}_{3^{3}}$. The orbits of this action are:  $\mathcal{O}_{3, 1}=\lbrace 1, 2, 4, \cdots, 25, 26 \rbrace$, $\mathcal{O}_{3, 3}=\lbrace 3, 6, 12, 15, 21, 24 \rbrace$, $\mathcal{O}_{3, 3^2}=\lbrace 9, 18 \rbrace$ and $\mathcal{O}_{3, 3^3}=\lbrace 0 \rbrace$. Therefore, by observation \eqref{4}, $$\Gamma(\mathbb{Z}_{3^3})_{\vee}=[\Gamma(\mathcal{O}_{3, 1}), \Gamma(\mathcal{O}_{3,3}), \Gamma(\mathcal{O}_{3,3^2}), \Gamma(\mathcal{O}_{3,3^3})],$$
where $\Gamma(\mathcal{O}_{3, 1})= \overline{K}_{18}$, $\Gamma(\mathcal{O}_{3, 3})= \overline{K}_{6}$, $\Gamma(\mathcal{O}_{3, 3^2})= {K}_{2}$, and $\Gamma(\mathcal{O}_{3, 3^3})= {K}_{1}$ as shown in the Figure \ref{b}.
\end{exm} 

\section{Threshold graphs realized by reduced, local and mixed rings}

The Artinian decomposition of $R =R_1\times \cdots \times R_r \times {F}_1 \times \cdots \times {F}_s$, where each $R_i$ is a finite commutative local ring with unity and ${F}_j$ is a field, $1 \leq i \leq r$ and $1 \leq j \leq s$. A ring $R$ is said to be product of local rings if $s =0$ and $r \geq 2$. A ring $R$ is said to be a \textit{mixed ring} if $r\geq1$ and $s\geq1$. A ring $R$ is said to be a \textit{reduced ring} if $r =0$ and $s \geq 1$.

This section is dedicated to study the realization of zero-divisor graphs as threshold graphs for reduced rings, product of local rings and mixed rings.

An alternating $4-cycle$ of a graph ${G} = (V, E)$ is a configuration consisting of four distinct vertices $a, b, c, d$ such that $(a,b), (c,d) \in E$  and $(a,c), (b,d) \notin E$. By considering the presence or absence of edges $(a,d)$ and $(b,c)$, we see that the vertices of an alternating $4$-cycle induce a path $P_4$, a cycle $C_4$, or a matching $2K_2$.
 
\begin{figure}[H]
  \centering
\includegraphics[scale=.3]{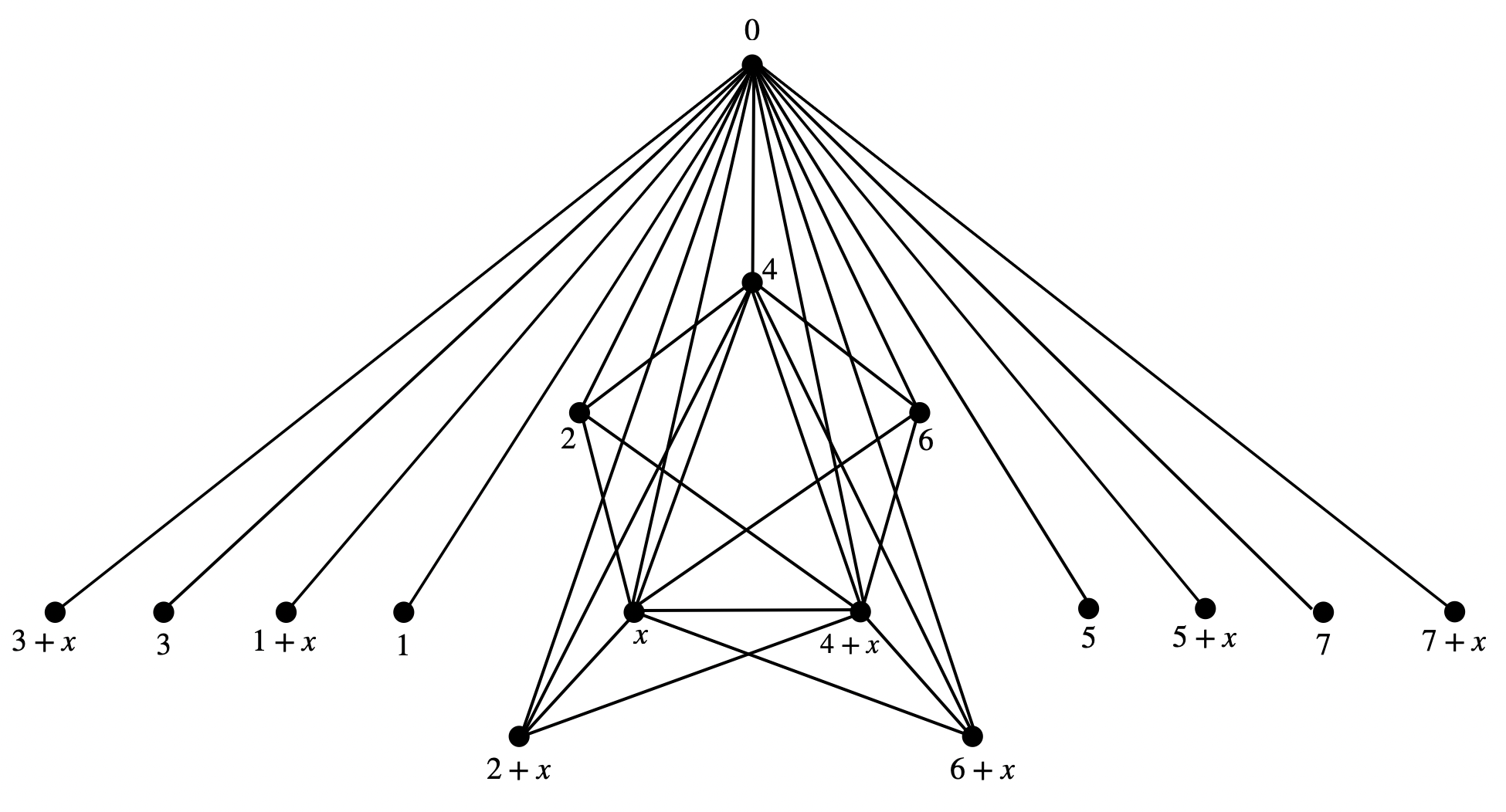}
\caption{$\Gamma\big(\mathbb{Z}_8[x]/(x^2 ,2x)\big)$}.
\label{c}
\end{figure}

Since threshold graphs has been characterised by many ways. One of the characterizations of a threshold graph [\cite{MP}, Theorem 1.2.4] is presented in the following result.
\begin{thm}
For a graph ${G} = (V, E ),$ the following are equivalent,\\
1. ${G}$ is a threshold graph,\\
2. ${G}$ does not have an alternating $4-cycle$.
\end{thm}
Let $R \cong {F}_{q_{1}}\times {F}_{q_{2}} \times \cdots \times {F}_{q_{n}}$ be a reduced ring. In the following characterization, we discuss the realization of $\Gamma(R)$ as a threshold graph.
  \begin{thm} 
  Let $R \cong \mathbb{F}_{q_{1}}\times {F}_{q_{2}}\times \cdots \times {F}_{q_{n}}$. Then $\Gamma(R)$ is a connected threshold graph if and only if $R\cong {F}_{q_{i}}$, or $R \cong  {F}_{q_{i}} \times {F}_{q_{j}}$, where ${F}_{q_{i}} \cong \mathbb{Z}_2$ for any $i, j \in \{1, 2, \cdots, n\}$.
  \label{11}
  \end{thm} 
  \begin{prf}
   Let $R$ be a ring such that $R\cong {F}_{q_i}$ for some $i \in \{1, 2, \cdots, n\}$. Then $\Gamma(R)$ is a complete bipartite graph isomorphic to $ K_{1,q_{i} - 1}$. This implies, $\Gamma(R)$ is a star graph. So, $\Gamma(R)$ is a connected threshold graph.
   
   Let $R$ be a ring such that $R \cong  {F}_{q_i} \times {F}_{q_j}$ for some  $i, j \in \{1, 2, \cdots, n\}$ with ${F}_{q_{i}} \cong \mathbb{Z}_2$. Then $|V(\Gamma(R))| = 2 {q_i}$.  Let $V(\Gamma(R)) = X\bigcup Y$, where $X$ represents the independent part of the graph whose vertices are of the form $(0, a_k)$ and $(1, a_k)$, $ k \in \{1, 2, \cdots, {q_j -1}\rbrace$. The set $Y$ represents complete part of the graph which consists of only two vertices  $(1,0)$ and $(0,0)$. 
    
    To prove $\Gamma(R)$ is threshold, it is suffices to show that $\Gamma(R)$ does not contain an alternating $4$-cycle, that is, the graph $\Gamma(R)$ has no induced subgraph of four vertices isomorphic to $P_4$, $C_4$ or $2K_2$.
    
    Let $a, b, c, d$ be any four distinct vertices of $\Gamma(R)$. We consider the following cases.\\
 {\bf Case - 1.} If $a, b, c, d \in X$, then they are all isolated vertices.
 This implies, $\Gamma(R)$ has no alternating $4$-cycle. Therefore, the result follows.\\
 {\bf Case - 2.}   If $a, b, c \in X$ and $d \in Y$, then the only possibilities are as follows,
 
 \begin{tikzpicture}[scale=0.8]
\draw[fill=black] (1,0) circle (3pt);
\draw[fill=black] (1,1) circle (3pt);
\draw[fill=black](2,0) circle (3pt);
\draw[fill=black](2,1) circle (3pt);
\node at (0.7,0) {c};
\node at (2.3,0) {d};
\node at (0.7,1) {a};
\node at (2.3,1) {b};
\draw[thick] (2,0) -- (1,0);
\draw[thick] (2,0) -- (2,1);
\draw[thick] (2,0) -- (1,1);
\end{tikzpicture}
\hfill
        \begin{tikzpicture}[scale=0.8]
\draw[fill=black] (1,0) circle (3pt);
\draw[fill=black] (1,1) circle (3pt);
\draw[fill=black] (2,0) circle (3pt);
\draw[fill=black] (2,1) circle (3pt);
\node at (0.7,0) {c};
\node at (2.3,0) {d};
\node at (0.7,1) {a};
\node at (2.3,1) {b};
\draw[thick] (2,0) (1,0);
\draw[thick] (2,0) -- (2,1);
\draw[thick] (2,0) -- (1,1);
\end{tikzpicture}
\hfill
        \begin{tikzpicture}[scale=0.8]
\draw[fill=black] (1,0) circle (3pt);
\draw[fill=black] (1,1) circle (3pt);
\draw[fill=black] (2,0) circle (3pt);
\draw[fill=black] (2,1) circle (3pt);
\node at (0.7,0) {c};
\node at (2.3,0) {d};
\node at (0.7,1) {a};
\node at (2.3,1) {b};
\draw[thick] (2,0) (1,0);
\draw[thick] (2,0) (2,1);
\draw[thick] (2,0) -- (1,1);
\end{tikzpicture}
\hfill
        \begin{tikzpicture}[scale=0.8]
\draw[fill=black] (1,0) circle (3pt);
\draw[fill=black] (1,1) circle (3pt);
\draw[fill=black] (2,0) circle (3pt);
\draw[fill=black] (2,1) circle (3pt);
\node at (0.7,0) {c};
\node at (2.3,0) {d};
\node at (0.7,1) {a};
\node at (2.3,1) {b};
\draw[thick] (2,0) (1,0);
\draw[thick] (2,0) (2,1);
\draw[thick] (2,0) (1,1); 
\end{tikzpicture}

but none of them is isomorphic to $P_4$, $C_4$ or $2K_2$. Therefore,  $\Gamma(R)$ is a connected threshold graph.\
    
{\bf Case - 3.} If $a, b \in X$ and $ c, d \in Y$, where $c = (1,0)$ and $ d = (0,0)$. By definition, the vertex $(0,0)$ is adjacent to three other vertices. Consequently, no two disjoint edges exists and therefore no alternating $4-cycle$ is possible. Thus, $\Gamma(R)$ is threshold.

Conversely, suppose that $\Gamma(R)$ is a connected threshold graph, we show that $R\cong {F}_{q_i}$, or $R \cong  {F}_{q_i} \times {F}_{q_j}$, where ${F}_{q_{i}} \cong \mathbb{Z}_2$ for any $i, j \in \{1, 2, \cdots, n\}$. We consider the following cases.\

{\bf Case - 1.} Suppose $R \cong  {F}_{q_1} \times {F}_{q_2}$,   $ 2<q_1 \leq q_2$.  Then, $|V(\Gamma(R))| = q_1 q_2$. Since $(a,b), (a,c), (c,d),\\ (b,d) \in E(\Gamma(R))$, where $ a = (1,0), b = (0, y), c = (0,1)$ and $d = (x,0)$, $x \in {F}_{q_1}\setminus \{0, 1\}$, $y \in {F}_{q_2}\setminus \{0, 1\}$. This implies, $\Gamma(R)$ has an induced subgraph of order $4$ isomorphic to $C_4$, a contradiction. Thus, $R  \ncong  {F}_{q_1} \times {F}_{q_2}$,   $ 2<q_1 \leq q_2$.\

{\bf Case - 2.} Suppose $ R \cong {F}_{q_{1}}\times {F}_{q_{2}}\times \cdots \times {F}_{q_{n}}$ with $q_1 \leq q_2 \leq \cdots \leq q_n$ and $n \geq 3$. Then, $\Gamma(R)$ always contains an induced subgraph of the form $P_4 = (1, 0, 1, 0, \cdots , 0) \sim (0, 1, 0, 0, \cdots,0) \sim (1, 0, 0, 0, \cdots,0) \sim (0, 1, 1, 0, \cdots,0)$, or $2K_2 = \{(1, 0, 1, 0, \cdots , 0) \sim (0, 1, 0, 1, \cdots, 0)$, $(0, 1, 1, 0, \cdots , 0)\\ \sim (1, 0, 0, 1, \cdots,0)\}$ or $C_4 = (a,b), (a,c), (c,d), (b,d)$, where $ a = (1,0,0, \cdots, 0), b = (0, y, 0, \cdots, 0), c = (0,1, 0, \cdots, 0)$ and $d = (x,0, 0, \cdots, 0)$, $x \in {F}_{q_1}\setminus \{0, 1\}$, $y \in {F}_{q_2}\setminus \{0, 1\}$, a contradiction. This implies, $ R \neq {F}_{q_{1}}\times {F}_{q_{2}}\times \cdots \times {F}_{q_{n}}$, for $q_1 \leq q_2 \leq \cdots \leq q_n$.

Therefore, $\Gamma(R)$ fails to be threshold for cases $1$ and $2$. We conclude that either  $R\cong {F}_{q_i}$ or $R \cong  {F}_{q_{i}} \times {F}_{q_j}$.
\end{prf} 
Now, we consider the case of local rings. We show that if $R$ is a local ring, then $\Gamma(R)$ is not always a threshold graph. 
Let $R \cong \mathbb{Z}_8[x]/(x^2 ,2x), \mathbb{Z}_2[x]/(x^3, 2x^2,2x)$ or $\mathbb{Z}_2[x,y]/(x^3,xy,y^2)$. Then the zero-divisor graph realized by $\mathbb{Z}_8[x]/(x^2 ,2x)$ is shown in the Figure \ref{c}. If we relabel the vertices of $\Gamma(\mathbb{Z}_8[x]/(x^2 ,2x))$ by elements of rings $\mathbb{Z}_2[x]/(x^3, 2x^2,2x)$ and $\mathbb{Z}_2[x,y]/(x^3,xy,y^2)$, then we obtain isomorphic copies of $\Gamma(\mathbb{Z}_8[x]/(x^2 ,2x))$ as zero-divisor graphs of $\mathbb{Z}_2[x]/(x^3, 2x^2,2x)$ and $\mathbb{Z}_2[x,y]/(x^3,xy,y^2)$.  It can be easily verified from the Figure \ref{c}, that $\Gamma\big(\mathbb{Z}_8[x]/(x^2 ,2x)\big)$ is a threshold graph.

\begin{figure}[H]
  \centering
\begin{tikzpicture}[scale=0.7]
\draw[fill=black] (0,4) circle (2pt);
\draw[fill=black] (-9 ,0) circle (2pt);
\draw[fill=black] (-7,0) circle (2pt);
\draw[fill=black] (-5,0) circle (2pt);
\draw[fill=black] (-4,0) circle (2pt);
\draw[fill=black] (-3,0) circle (2pt);
\draw[fill=black] (-2,-2) circle (2pt);
\draw[fill=black] (0,0) circle (2pt);
\draw[fill=black] (-1.5,1.5) circle (2pt);
\draw[fill=black] (1.5,1.5) circle (2pt);
\draw[fill=black] (3,0) circle (2pt);
\draw[fill=black] (2,-2) circle (2pt);
\draw[fill=black] (4,0) circle (2pt);
\draw[fill=black] (5,0) circle (2pt);
\draw[fill=black] (7,0) circle (2pt);
\draw[fill=black] (9,0) circle (2pt);
\node at (0,4.4) {0};
\node at (-9,-0.4) {1+3x};
\node at (-7,-0.4) {1+2x};
\node at (-5,-0.4) {1+x};
\node at (-4,-0.4) {1};
\node at (-3,-0.4) {2+x};
\node at (-2,-2.4) {2+3x};
\node at (0,-0.4) {2x};
\node at (-1.5,1.1) {x};
\node at (1.5,1.1) {3x};
\node at (3,-0.4) {2+2x};
\node at (2,-2.4) {2};
\node at (4,-0.4) {3};
\node at (5,-0.4) {3+x};
\node at (7,-0.4) {3+2x};
\node at (9,-0.4) {3+3x};

\draw[thick] (0,0) -- (-1.5,1.5);
\draw[thick] (0,0) -- (1.5,1.5);
\draw[thick] (0,0) -- (-3,0);
\draw[thick] (0,0) -- (3,0);
\draw[thick] (0,0) -- (-2,-2); 
\draw[thick] (0,0) -- (2,-2);
\draw[thick] (-1.5,1.5) -- (1.5,1.5);
\draw[thick] (-3,0) -- (-2,-2); 
\draw[thick] (3,0) -- (2,-2);
\draw[thick] (0,4) -- (-9,0);
\draw[thick] (0,4) -- (-7,0); 
\draw[thick] (0,4) -- (-5,0);
\draw[thick] (0,4) -- (-4,0);
\draw[thick] (0,4) -- (-3,0); 
\draw[thick] (0,4) -- (-2,-2);
\draw[thick] (0,4) -- (0,0);
\draw[thick] (0,4) -- (-1.5,1.5);
\draw[thick] (0,4) -- (1.5,1.5);
\draw[thick] (0,4) -- (3,0);
\draw[thick] (0,4) -- (4,0);
\draw[thick] (0,4) -- (2,-2);
\draw[thick] (0,4) -- (5,0);
\draw[thick] (0,4) -- (7,0);
\draw[thick] (0,4) -- (9,0);
\end{tikzpicture}
\caption{$\Gamma\big(\mathbb{Z}_4[x]/(x^2)\big)$}
\label{d}
\end{figure}
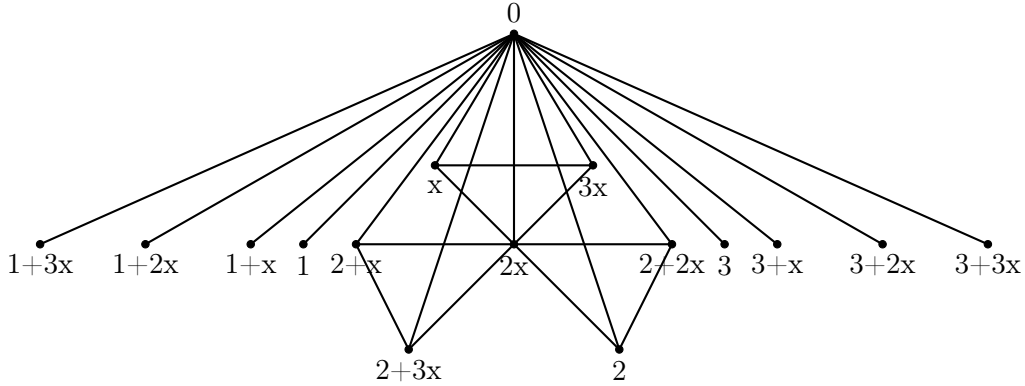

On the other hand, if $R \cong \mathbb{Z}_4[x]/(x^2), \mathbb{Z}_2[x,y]/(x^2,y^2)$ or $\mathbb{Z}_4[x,y]/(x^2, y^2, xy-2,2x, 2y)$, then the graph realized by these rings are isomorphic to each other. The graph realized by the ring $\mathbb{Z}_4[x]/(x^2)\big)$ is shown in the Figure \ref{d}. The graph $\Gamma\big(\mathbb{Z}_4[x]/(x^2)\big)$ fails to be threshold, since it contains an induced subgraph isomorphic to $2K_2$, that is the edges $(x,3x)$ and $(2,2+2x)$ are disjoint in $\Gamma\big(\mathbb{Z}_4[x]/(x^2)\big)$.
Thus from preceding examples it is hard to characterize local rings such that their zero-divisor graphs are threshold. Below, we present some classes of local rings of which associated zero-divisor graphs are threshold. 
\begin{thm}
For each prime $p$ and positive integer $\alpha$, the zero-divisor graph realized by the ring $\mathbb{Z}_{p^{\alpha}}[x]/(x^2,px)$ is a connected threshold graph.
\label{14}
\end{thm}
\begin{prf} Let $R \cong \mathbb{Z}_{p^{\alpha}}[x]/(x^2,px)$. Then orbits of the group action $Aut(\Gamma(R)) \times R \longrightarrow R$ are as follows,

$\mathcal{O}_{\alpha,p^i} = \{p^ia+bx :0\neq a\in \mathbb{Z}_{p^{\alpha}}, ~ b\in \mathbb{Z}_p ~ \mbox {and} ~ (a,p)=1\}$, where $0\leq i\leq \alpha -2$, and

    $\mathcal{O}_{\alpha,p^{\alpha-1}} = \{p^{\alpha-1}a+bx\neq 0: ~a\in \mathbb{Z}_{p^{\alpha}}, ~ b\in \mathbb{Z}_p ~ \mbox {and} ~ (a,p)=1\}.$
    
\noindent Moreover, for $0\leq i\leq \alpha-2$,  
$|\mathcal{O}_{\alpha,p^i}| = \phi\big(\frac{p^{\alpha+1}}{p^i}\big)$ and $|\mathcal{O}_{\alpha,p^{\alpha-1}}| = \phi\big(\frac{p^{\alpha+1}}{p^{\alpha-1}}\big)+p-1=p^2-1$. Note that $\mathcal{O}_{\alpha,p^\alpha}=\lbrace0\rbrace$. The vertex set $V\big(\Gamma(R)\big)= \big(\bigcup\limits_{i=0}^{\alpha} \mathcal{O}_{\alpha,p^{i}}\big)$. To prove the desired result, we consider the following  cases.\\
\textbf{Case - 1.} $\alpha=2k$ for some $k \geq 1$. Then the subgraphs associated with orbits $X =\big(\bigcup\limits_{i=0}^{k-1} \mathcal{O}_{\alpha,p^{k+i}}\big)$ represents complete part of the graph $\Gamma\big(R\big)$ and the subgraphs associated with orbits
  $Y = \big(\bigcup\limits_{i=0}^{k-1} \mathcal{O}_{\alpha,p^{i}}\big)$ represents the independent part of $\Gamma(R)$ and each element of $X$ and $Y$ are connected to $ \mathcal{O}_{\alpha,p^{\alpha}}=\lbrace0\rbrace$. Furthermore, if $r\in \mathcal{O}_{\alpha,p^{k+i}}$ and $s\in  \mathcal{O}_{\alpha, p^{i}}$ for some $i$, then $r$ and $s$ are adjacent in the graph if and only if $p^\alpha$ divides $p^{2i + k}$. This implies,
  
$|V\big(\Gamma(R)\big)| =  \sum\limits_{i=0}^{k-1} |\mathcal{O}_{\alpha,p^{i}}| + \sum\limits_{i=0}^{k-2}|\mathcal{O}_{\alpha,p^{k+i}}| +\big(|\mathcal{O}_{\alpha,p^{2k-1}}|+(p-1)\big)+ 1 = p^{\alpha+1}.$\\
 To prove that $\Gamma(R)$ is a threshold graph it suffices to show that $\Gamma(R)$ has no induced subgraph of the form $C_4$, $P_4$ or $2K_2$, that is, it has no alternating $4$-$cycle$.
 
 Let $a, b, c, d$ be four distinct vertices of $\Gamma(R)$.\\
 \textbf{Subcase - 1.} If $a = 0$ and $b, c ,d \in X \cup Y$, then $(b,0), (c,0), (d,0) \in E\big(\Gamma(R)\big)$. Thus no two disjoint edges are present and therefore no alternating $4-cycle$ is possible.\\
 \textbf{Subcase - 2.} If all $a, b, c ,d \in X \cup Y$, that is, no one among $a, b, c ,d$ represents the vertex $0$, then we have the following cases .\\
 \textbf{Subcase - 2.1.} Suppose that all four vertices $a, b, c ,d \in X$. Then they are all adjacent to each other and hence no alternating $4-cycle$ is possible.\\
 \textbf{Subcase - 2.2.} Suppose that all four vertices $a, b, c ,d \in Y$. Then these vertices are all isolated and therefore no alternating $4-cycle$  is possible.\\
 \textbf{Subcase - 2.3.} Suppose $a \in X$ and $b ,c, d \in Y$. Then in $\Gamma(R)$, there may not be any edge from $a$ to $b, c, d$. Without of loss of generality, the other possibilities are either $(a, d) \in E(\Gamma(R))$ or $(a, b), (a, d) \in E(\Gamma(R))$ or $(a, b), (a, c), (a, d) \in E(\Gamma(R))$ as follows.
 
 \begin{tikzpicture}[scale=0.8]
\draw[fill=black] (1,0) circle (3pt);
\draw[fill=black] (1,1) circle (3pt);
\draw[fill=black] (2,0) circle (3pt);
\draw[fill=black] (2,1) circle (3pt);
\node at (0.7,0) {c};
\node at (2.3,0) {d};
\node at (0.7,1) {a};
\node at (2.3,1) {b};
\draw[thick] (2,0) (1,0);
\draw[thick] (2,0) (2,1);
\draw[thick] (2,0) (1,1); 
\end{tikzpicture}
\hfill
  \begin{tikzpicture}[scale=0.8]
\draw[fill=black] (1,0) circle (3pt);
\draw[fill=black] (1,1) circle (3pt);
\draw[fill=black] (2,0) circle (3pt);
\draw[fill=black] (2,1) circle (3pt);
\node at (0.7,0) {c};
\node at (2.3,0) {d};
\node at (0.7,1) {a};
\node at (2.3,1) {b};
\draw[thick] (2,0) (1,0);
\draw[thick] (2,0) (2,1);
\draw[thick] (2,0) -- (1,1);
\end{tikzpicture}
\hfill
        \begin{tikzpicture}[scale=0.8]
\draw[fill=black] (1,0) circle (3pt);
\draw[fill=black] (1,1) circle (3pt);
\draw[fill=black] (2,0) circle (3pt);
\draw[fill=black] (2,1) circle (3pt);
\node at (0.7,0) {c};
\node at (2.3,0) {d};
\node at (0.7,1) {a};
\node at (2.3,1) {b};
\draw[thick] (2,0) (1,0);
\draw[thick] (1,1) -- (2,0);
\draw[thick] (1,1) -- (2,1);
\end{tikzpicture}
\hfill
 \begin{tikzpicture}[scale=0.8]
\draw[fill=black] (1,0) circle (3pt);
\draw[fill=black] (1,1) circle (3pt);
\draw[fill=black](2,0) circle (3pt);
\draw[fill=black](2,1) circle (3pt);
\node at (0.7,0) {c};
\node at (2.3,0) {d};
\node at (0.7,1) {a};
\node at (2.3,1) {b};
\draw[thick] (1,1) -- (1,0);
\draw[thick] (1,1) -- (2,1);
\draw[thick] (1,1) -- (2,0);
\end{tikzpicture}\\
Thus no one among the above graphs is isomorphic to $C_4$, $P_4$ or $2K_2$. Therefore the result follows.\\
\textbf{Subcase - 2.4.} Suppose $a, b \in X$ and $c,d \in Y$, then the only possibilities are as follows. \
 
 \begin{tikzpicture}[scale=0.8]
\draw[fill=black] (1,0) circle (3pt);
\draw[fill=black] (1,1) circle (3pt);
\draw[fill=black] (2,0) circle (3pt);
\draw[fill=black] (2,1) circle (3pt);
\node at (0.7,0) {c};
\node at (2.3,0) {d};
\node at (0.7,1) {a};
\node at (2.3,1) {b};
\draw[thick] (2,0) (1,0);
\draw[thick] (2,0) (2,1);
\draw[thick] (2,1) -- (1,1); 
\end{tikzpicture}
\hfill
        \begin{tikzpicture}[scale=0.8]
\draw[fill=black] (1,0) circle (3pt);
\draw[fill=black] (1,1) circle (3pt);
\draw[fill=black] (2,0) circle (3pt);
\draw[fill=black] (2,1) circle (3pt);
\node at (0.7,0) {c};
\node at (2.3,0) {d};
\node at (0.7,1) {a};
\node at (2.3,1) {b};
\draw[thick] (2,1) -- (1,0);
\draw[thick] (1,1) (2,0);
\draw[thick] (1,1) -- (2,1);
\end{tikzpicture}
\hfill
 \begin{tikzpicture}[scale=0.8]
\draw[fill=black] (1,0) circle (3pt);
\draw[fill=black] (1,1) circle (3pt);
\draw[fill=black](2,0) circle (3pt);
\draw[fill=black](2,1) circle (3pt);
\node at (0.7,0) {c};
\node at (2.3,0) {d};
\node at (0.7,1) {a};
\node at (2.3,1) {b};
\draw[thick] (1,1) -- (1,0);
\draw[thick] (1,1) -- (2,1);
\draw[thick] (2,1) (2,0);
\draw[thick] (1,0) -- (2,1);
\end{tikzpicture}
\hfill
  \begin{tikzpicture}[scale=0.8]
\draw[fill=black] (1,0) circle (3pt);
\draw[fill=black] (1,1) circle (3pt);
\draw[fill=black] (2,0) circle (3pt);
\draw[fill=black] (2,1) circle (3pt);
\node at (0.7,0) {c};
\node at (2.3,0) {d};
\node at (0.7,1) {a};
\node at (2.3,1) {b};
\draw[thick] (2,1) -- (1,0);
\draw[thick] (2,1) -- (2,0);
\draw[thick] (2,1) -- (1,1);
\end{tikzpicture}
\hfill
 \begin{tikzpicture}[scale=0.8]
\draw[fill=black] (1,0) circle (3pt);
\draw[fill=black] (1,1) circle (3pt);
\draw[fill=black](2,0) circle (3pt);
\draw[fill=black](2,1) circle (3pt);
\node at (0.7,0) {c};
\node at (2.3,0) {d};
\node at (0.7,1) {a};
\node at (2.3,1) {b};
\draw[thick] (1,1) -- (1,0);
\draw[thick] (1,1) -- (2,1);
\draw[thick] (2,1) -- (2,0);
\draw[thick] (1,0) -- (2,1);
\end{tikzpicture}
\hfill
 \begin{tikzpicture}[scale=0.8]
\draw[fill=black] (1,0) circle (3pt);
\draw[fill=black] (1,1) circle (3pt);
\draw[fill=black](2,0) circle (3pt);
\draw[fill=black](2,1) circle (3pt);
\node at (0.7,0) {c};
\node at (2.3,0) {d};
\node at (0.7,1) {a};
\node at (2.3,1) {b};
\draw[thick] (1,1) -- (1,0);
\draw[thick] (1,1) -- (2,1);
\draw[thick] (1,1) -- (2,0);
\draw[thick] (1,0) -- (2,1);
\draw[thick] (2,0) -- (2,1);
\end{tikzpicture}\\
\textbf{Subcase - 2.5.} Suppose $a, b, c \in X$ and $d \in Y$, then the only possibilities are as follows.\

\begin{tikzpicture}[scale=0.8]
\draw[fill=black] (1,0) circle (3pt);
\draw[fill=black] (1,1) circle (3pt);
\draw[fill=black] (2,0) circle (3pt);
\draw[fill=black] (2,1) circle (3pt);
\node at (0.7,0) {c};
\node at (2.3,0) {d};
\node at (0.7,1) {a};
\node at (2.3,1) {b};
\draw[thick] (1,1) -- (1,0);
\draw[thick] (2,1) -- (1,0);
\draw[thick] (2,1) -- (1,1); 
\end{tikzpicture}
\hfill
        \begin{tikzpicture}[scale=0.8]
\draw[fill=black] (1,0) circle (3pt);
\draw[fill=black] (1,1) circle (3pt);
\draw[fill=black] (2,0) circle (3pt);
\draw[fill=black] (2,1) circle (3pt);
\node at (0.7,0) {c};
\node at (2.3,0) {d};
\node at (0.7,1) {a};
\node at (2.3,1) {b};
\draw[thick] (2,1) -- (1,0);
\draw[thick] (1,1) -- (1,0);
\draw[thick] (1,1) -- (2,1);
\draw[thick] (1,0) -- (2,0);
\end{tikzpicture}
\hfill
 \begin{tikzpicture}[scale=0.8]
\draw[fill=black] (1,0) circle (3pt);
\draw[fill=black] (1,1) circle (3pt);
\draw[fill=black](2,0) circle (3pt);
\draw[fill=black](2,1) circle (3pt);
\node at (0.7,0) {c};
\node at (2.3,0) {d};
\node at (0.7,1) {a};
\node at (2.3,1) {b};
\draw[thick] (1,1) -- (1,0);
\draw[thick] (1,1) -- (2,1);
\draw[thick] (2,1) -- (2,0);
\draw[thick] (1,0) -- (2,1);
\draw[thick] (1,0) -- (2,0);
\end{tikzpicture}
\hfill
  \begin{tikzpicture}[scale=0.8]
\draw[fill=black] (1,0) circle (3pt);
\draw[fill=black] (1,1) circle (3pt);
\draw[fill=black] (2,0) circle (3pt);
\draw[fill=black] (2,1) circle (3pt);
\node at (0.7,0) {c};
\node at (2.3,0) {d};
\node at (0.7,1) {a};
\node at (2.3,1) {b};
\draw[thick] (1,1) -- (1,0);
\draw[thick] (1,1) -- (2,1);
\draw[thick] (2,1) -- (2,0);
\draw[thick] (1,0) -- (2,1);
\draw[thick] (1,0) -- (2,0);
\draw[thick] (1,1) -- (2,0);
\end{tikzpicture}\\
Thus it follows that the graph $\Gamma(R)$ has no alternating $4-cycle$. Therefore, $\Gamma(R)$ is a connected threshold graph.\\
\textbf{Case - 2.}
$\alpha=2k-1$ for some $k \geq 1$. Then the subgraphs associated with orbits $X =\big(\bigcup\limits_{i=0}^{k-2} \mathcal{O}_{\alpha,p^{k+i}}\big)$ represents complete part of the graph $\Gamma\big(R\big)$ and the subgraphs associated with orbits
  $Y = \big(\bigcup\limits_{i=0}^{k-1} \mathcal{O}_{\alpha,p^{i}}\big)$ represents the independent part of $\Gamma(R)$ and each element of $X$ and $Y$ are connected to $ \mathcal{O}_{\alpha,p^{\alpha}}=\lbrace0\rbrace$. This implies,\
 
$|V\big(\Gamma(R)\big)| =  \sum\limits_{i=0}^{k-1} |\mathcal{O}_{\alpha,p^{i}}| + \sum\limits_{i=0}^{k-3}|\mathcal{O}_{\alpha,p^{k+i}}| +\big(|\mathcal{O}_{\alpha,p^{2k-2}}|+(p-1)\big)+ 1 = p^{\alpha+1}.$

We can prove that $ \Gamma\big(R\big)$ does not have an alternating $4$-cycle by considering four distinct vertices $a,b,c,d$ and using the same argument as above for the case $\alpha = 2k$.
\end{prf}
\begin{thm}
For each prime $p$, zero-divisor graphs realized by rings $\mathbb{Z}_{p}[x]/(x^p)$, $\mathbb{Z}_{p}[x, y]/(x^3, xy, y^2)$, and $\mathbb{Z}_{p^2}[x]/(px, x^2-p)$ are connected threshold graphs.
\label{ss}
\begin{prf} Let $R \cong \mathbb{Z}_{p}[x]/(x^p)$. Then orbits of the group action $Aut(\Gamma(R)) \times R \longrightarrow R$ are as follows,
$$\mathcal{O}_{p,x^i} = \{x^i(a_0+a_1x+\cdots+a_{p-(1+i)}x^{p-(1+i)}) : a_i\in \mathbb{Z}_{p},~ \mbox {and} ~ a_0\neq 0\},$$ where $0\leq i \leq p-1$. For $0\leq i \leq p-1$, $|\mathcal{O}_{p,x^i}| = \phi\big(\frac{p^{p}}{p^i}\big)$ and $|\mathcal{O}_{p,x^i}| = \lbrace0\rbrace$. The vertex set $V\big(\Gamma(R)\big)=X\bigcup Y$, where $X= \big(\bigcup\limits_{i=0}^{\frac{p-1}{2}}\mathcal{O}_{p,x^{i+\frac{p+1}{2}}}\big)$ and $Y= \big(\bigcup\limits_{i=0}^{\frac{p-1}{2}}\mathcal{O}_{p,x^{i}}\big)$, represents complete and independent parts of the graph $\Gamma(R)$, and $\mathcal{O}_{p,x^{p}}=\lbrace 0 \rbrace$. Also, two distinct vertices $r\in \mathcal{O}_{p,x^{i}}$ and $s\in \mathcal{O}_{p,x^{j} }$ for $i,~j\in\lbrace 0, 1, \cdots, \frac{p-1}{2}\rbrace$ are adjacent in $\Gamma(R)$ if and only if $i+j\geq p$.  Therefore, $$|V\big(\Gamma(R)\big)| = \sum\limits_{i=0}^{\frac{p-1}{2}}|\mathcal{O}_{p,x^{i+\frac{p+1}{2}}}|+ \sum\limits_{i=0}^{\frac{p-1}{2}}|\mathcal{O}_{p,x^{i}}|=p^p.$$\\
Moreover, for $R\cong \mathbb{Z}_{p}[x, y]/(x^3, xy, y^2)$ and for each $a_0,~a_1,~a_2,~b_1\in \mathbb{Z}_{p}$, the orbits of the group action $Aut(\Gamma(R)) \times R \longrightarrow R$ are as follows,
\begin{eqnarray*}
\mathcal{O}_{p, x} &=& \{a_{1}x+a_{2}x^2+b_{1}y : ~~a_1\neq 0\},\\
\mathcal{O}_{p, y} &=& \{a_{2}x^2+b_{1}y\neq 0 \},\\
\mathcal{O}_{p, 1} &=& \{a_{0}+a_{1}x+a_{2}x^2+b_{1}y: ~~a_{0}\neq 0\},\\
\mathcal{O}_{p, p} &=& \{0\},
\end{eqnarray*}
where $x, y, p$ and $1$ are representatives of orbits $\mathcal{O}_{p, x}$, $\mathcal{O}_{p, y}$, $\mathcal{O}_{p, 1}$ and $\mathcal{O}_{p, p}$ with  $|\mathcal{O}_{p, x}|=p^3-p^2$, $|\mathcal{O}_{p, y}|=p^2-1$, $|\mathcal{O}_{p, 1}|=p^4-p^3$ and $|\mathcal{O}_{p, p}|=1$. The vertex set $V\big(\Gamma(R)\big)=X\bigcup Y$,  where $X= \big(\mathcal{O}_{p, y}\big)$ and $Y= \big(\mathcal{O}_{p, x}\cup \mathcal{O}_{p, 1}\big)$ represents complete and independent parts of the graph $\Gamma(R)$, and $\mathcal{O}_{p, p}=\lbrace 0 \rbrace$. Furthermore, every vertex in $\mathcal{O}_{p, y}$ is adjacent to each vertex in $\mathcal{O}_{p, x}$ and the vertex $0$ is adjacent to each vertex of $\Gamma(R)$. Therefore, $|V\big(\Gamma(R)\big)| = p^4$.\

Finally, for $R\cong \mathbb{Z}_{p^2}[x]/(px, x^2-p)$, the orbits of the group action $Aut(\Gamma(R)) \times R \longrightarrow R$ are as follows,
\begin{eqnarray*}
\mathcal{O}_{p, x} &=& \{pa+bx : ~0\neq b\in\mathbb{Z}_{p},~a\in\mathbb{Z}_{p}\},\\
\mathcal{O}_{p, p} &=& \{pa : ~0\neq a\in\mathbb{Z}_{p}\},\\
\mathcal{O}_{p, 1} &=& \{a+bx: ~0\neq a\in\mathbb{Z}_{p^2},~b\in\mathbb{Z}_{p}~and~ (a, p)=1\}\\
\mathcal{O}_{p, p^2} &=& \{0\},
\end{eqnarray*}
where $|\mathcal{O}_{p, x}|=p^2-p$, $|\mathcal{O}_{p, p}|=p-1$, $|\mathcal{O}_{p, 1}|=p^3-p^2$ and $|\mathcal{O}_{p, p^2}|=1$. The vertex set $V\big(\Gamma(R)\big)=X\bigcup Y$, where $X= \big(\mathcal{O}_{p, p}\big)$ and $Y= \big(\mathcal{O}_{p, x}\cup \mathcal{O}_{p, 1}\big)$  represents complete and independent parts of the graph $\Gamma(R)$ and $\mathcal{O}_{p, p^2}=\lbrace 0 \rbrace$. From the structure of the graph, every vertex in $\mathcal{O}_{p, 1}$ is adjacent to the vertex $0$ only, whereas every vertex in $\mathcal{O}_{p, p}$ is adjacent to each vertex of $\mathcal{O}_{p, x}$ and to the vertex $0$. Thus, $|V\big(\Gamma(R)\big)| = p^3$.\

By taking any four distinct vertices a, b, c and d from vertex set of each of the graphs realized by rings $(\mathbb{Z}_{p}[x]/(x^p))$, $\mathbb{Z}_{p}[x, y]/(x^3, xy, y^2)$, $\mathbb{Z}_{p^2}[x]/(px, x^2-p)$ and use the similar approach as in Theorem \eqref{14}, it can be shown that zero divisor graphs realized by above mentioned rings does not contain an alternating $4-cycle$ and therefore the result follows.
\end{prf}
\end{thm}
\begin{cor}
For each prime $p$ and positive integer $\alpha$, the zero-divisor graph realized by a local ring $\Z_{p^{\alpha}}$ is a connected threshold graph.
\end{cor}

\textbf{Problem.}  Characterise all finite local rings such that $\Gamma(R)$ is a threshold graph.    

In the following result, we consider the Artinian decomposition of $R$ for $r \geq 2$ and $s = 0$, that is, we investigate the structure of $\Gamma(R)$ when $R$ is a product of local rings only. 
\begin{lem}
 If $R\cong{R}_1 \times {R}_2 \times \cdots \times {R}_r$, where $r\geq2 $ and each ${R}_i, 1\leq i \leq r$ is a finite local ring (not a field), then $\Gamma(R)$ is not a threshold graph.
 \label{12}
\end{lem} 
\begin{prf}
 Since $R$ is a product of at least two finite local rings, therefore $\Gamma(R)$ always has an induced subgraph either of the form $P_4$ or $C_4$ or $2K_2$ described as follows.

 \centering
 
 \begin{tikzpicture}[scale=0.8]
\draw[fill=black] (0,0) circle (3pt);
\draw[fill=black] (4,0) circle (3pt);
\draw[fill=black] (8,0) circle (3pt);
\draw[fill=black] (12,0) circle (3pt);
\node at (0,-0.4) {(1, ${b}_1$, 0, $\cdots$ ,0)};
\node at (4,-0.4) {(0, ${b}_2$, 0, $\cdots$ ,0)};
\node at (8,-0.4) {(${a}_1$, 0, 0, $\cdots$ ,0)};
\node at (12,-0.4) {(${a}_2$, 1, 0, $\cdots$ ,0)};
\draw[thick] (0,0) -- (4,0) -- (8,0) -- (12,0);
\end{tikzpicture}
\hfill
 \begin{tikzpicture}[scale=0.8]
\draw[fill=black] (0,0) circle (3pt);
\draw[fill=black] (4,0) circle (3pt);
\draw[fill=black] (0,4) circle (3pt);
\draw[fill=black] (4,4) circle (3pt);
\node at (0,-0.4) {(1, 0, 0, $\cdots$ ,0)};
\node at (4,-0.4) {(0, 1, 0, $\cdots$ ,0)};
\node at (0,4.4) {(0, ${b}_1$, 0, $\cdots$ ,0)};
\node at (4,4.4) {(${a}_1$, 0, 0, $\cdots$ ,0)};
\draw[thick] (0,0) -- (4,0) -- (4,4) -- (0,4);
\draw[thick] (0,0) -- (0,4);
\end{tikzpicture}
 \hfill
 \begin{tikzpicture}[scale=0.8]
\draw[fill=black] (0,0) circle (3pt);
\draw[fill=black] (4,0) circle (3pt);
\draw[fill=black] (0,4) circle (3pt);
\draw[fill=black] (4,4) circle (3pt);
\node at (0,-0.4) {(1, 0, 0, $\cdots$ ,0)};
\node at (0,4.4) {(0, 1, 0, $\cdots$ ,0)};
\node at (4,-0.4) {(${a}_1$, ${b}_1$, 0, $\cdots$ ,0)};
\node at (4,4.4) {(${a}_2$, ${b}_2$, 0, $\cdots$ ,0)};
\draw[thick] (0,0) -- (0,4);
\draw[thick] (4,0) -- (4,4);
\end{tikzpicture}\

Note that $a_i \in Z^*(R_1)$ and $b_i \in Z^*(R_2)$. where $Z^*(R_1)$ and $Z^*(R_2)$ denotes the set of non-zero zero-divisors of rings $R_1$ and $R_2$.
\end{prf}
In the following result, we classify all finite commutative rings $R$ such that $\Gamma(R)$ is not a threshold graph.
\begin{thm}
 If $R \cong {R}_1 \times {R}_2 \times \cdots \times {R}_r \times {{F}_1} \times {{F}_2} \times \cdots \times {{F}_s}$, where $r$ and $s$ are positive integers, each $R_i$ is a finite commutative local ring (not a field), and each ${F}_i$ is a finite field, then $\Gamma(R)$ is not threshold.
 \end{thm}
 \begin{prf}
 The proof follows from Theorem \eqref{11} and Lemma \eqref{12}.
 \end{prf}
                                                                                                                                                                                                                                                                                                                                                                                                                                                                                                                                                                                                                                                                                  
                                                                                                                                                                                                                                                                                                                                                                                                                                                                                                                                                                                                                                                                                  {\bf Conclusion:} In this research article, we presented some new classes of threshold graphs. As discussed above, the construction of new classes of threshold graphs are effective in a sense that it assists to advance in coding theory and spectral analysis. We obtained the results related to reduced and local rings. In fact, we characterized all reduced rings of which zero-divisor graphs are realized as threshold graphs. Furthermore, we obtained classes of threshold graphs realized by some classes of local rings and provided an example which illustrates that a zero-divisor graph associated with some local ring is not always a threshold graph. So, it is interesting to determine all classes of local rings of which zero-divisor graphs are always threshold. This investigation has been included in the paper as a problem for future research work.

                                                                                                                                                                                                                                                                                                                                                                                                                                                                                                                                                                                                                                                                                  {\bf Acknowledgement}\\ 
                                                                                                                                                                                                                                                                                                                                                                                                                                                                                                                                                                                                                                                                                  We thank the anonymous reviewers for their careful reading of the manuscript and their many insightful comments and suggestions. The second author's research is supported by the University Grants Commission, Govt. of India under UGC-Ref. No. 1337/(CSIR-UGC NET JUNE 2019).

                                                                                                                                                                                                                                                                                                                                                                                                                                                                                                                                                                                                                                                                                  \textbf{Declaration of competing interest}\\
                                                                                                                                                                                                                                                                                                                                                                                                                                                                                                                                                                                                                                                                                  There is no conflict of interest to declare.

\textbf{Data Availability.}\\
Data sharing not applicable to this article as no datasets were generated or analysed during the current study.

\end{document}